\documentclass[12pt]{amsart}
\usepackage{amssymb}

  \theoremstyle{definition}
  \newtheorem{theorem}{Theorem}[section] 
  \newtheorem{corollary}[theorem]{Corollary} 
  
  \newtheorem{lemma}[theorem]{Lemma}
  \newtheorem{proposition}[theorem]{Proposition}
  
  \newtheorem{remark}[theorem]{Remark}
  \newtheorem{example}[theorem]{Example}
  
  \newtheorem{definition}[theorem]{Definition}
  \newtheorem{question}[theorem]{Question}
  \theoremstyle{remark}
  \numberwithin{equation}{section}


\begin{document}

\title[Complemented Subspaces of $L_p$ Determined by Partitions and Weights]
   {Subspaces of $L_p$, $p>2$, Determined by Partitions and Weights}
\vskip 2 cm
\author {Dale E. Alspach}
\address{Department of Mathematics, Oklahoma State
   University, Stillwater, OK 74078}
\email{alspach@math.okstate.edu}
\author {Simei Tong}
\address{Department of Mathematics, University of Wisconsin-Eau
   Claire, Eau Claire, WI 54701}
\email{tongs@uwec.edu}
\subjclass[2000]{Primary 46B20 Secondary 46E30}

\today

\begin{abstract}Many of the
known complemented subspaces of $L_p$ have realizations
as sequence spaces. In this paper a systematic approach to defining
these spaces which uses partitions and weights is introduced. This approach
gives a
unified description of many well-known complemented subspaces of
$L_p$. It is proved that the class of spaces with such norms is stable under
$(p,2)$ sums.
 By introducing the notion of an envelope norm, we obtain a necessary condition 
for a Banach sequence space with norm given by partitions and weights
 to be isomorphic to a subspace of $L_p$. Using this we define  
 a space $Y_n$ with norm given by partitions and weights 
with distance to any subspace of $L_p$ growing with $n$. This allows
us to construct an example of a Banach space with norm
given by partitions and weights which is not isomorphic to
 a subspace of $L_p$.
\end{abstract}
\maketitle
 

\section{Introduction} Prior to Rosenthal's 1970 paper, {\bf [R]}, only a few
complemented subspaces of $L_p$ were known, $\ell_p$, $\ell_2$, 
$\ell_p\oplus \ell_2$, $(\sum \ell_2)_p$ and $L_p$ itself, {\bf [P]},
{\bf [F]}. Rosenthal's
paper added several new spaces but more importantly it was seminal. In
1975 Schechtman, {\bf [S]}, combined Rosenthal's results with a tensor product 
construction to show that there were infinitely many isomorphically
distinct complemented subspaces of $L_p$. A few years later, it was
shown by Bourgain, Rosenthal 
and Schechtman that up to isomorphism, there are uncountably
 many complemented subspaces of $L_p$ {\bf{[B-S-R]}}. Recently
the first author proposed a new approach to describe the
 complemented subspaces of $L_p[0,1], p>2$. Define for any partition 
$P=\{ N_i\}$ of $\mathbb{N}$ and weight function $ W : {\mathbb{N}}
 \rightarrow { ( 0, 1 ]}$
  $$\begin{displaystyle} \|(a_i)\|_{P,W} = \left(\sum_{i} \left (
 \sum_{j\in {N_i}} {a_j}^2 w_{j}^{2}\right)^{\frac{p}{2}}\right)
^{\frac{1}{p}} \end{displaystyle}$$
Now suppose that ${\mathcal{P}}=(P_k, W_k)_{k\in K}$ is a family of pairs of
partitions and functions as above.
Define \begin{eqnarray}
 \|(a_i)\|_{\mathcal{P}} =  \sup_{k\in K}
 \|(a_i)\|_{(P_k, W_k)} 
\end{eqnarray}\label{norm-0}
There are two fundamental questions 
which we begin to study in this paper. What conditions on 
$(P_k, W_k)_{k \in K}$ imply that $\|(a_i)\|_{\mathcal{P}}$  
 defines a norm on a 
 space of sequences $X$ so that $X$ is isomorphic to a complemented
 subspace of $L_p[0, 1]$? Is every complemented subspace of $L_p$ 
(other than $L_p$) isomorphic to a space of this form?


This paper includes four major sections in addition to this introduction.
Unless otherwise noted we will assume that $p>2$ throughout. We will
also assume that the scalar field is $\mathbb{R}$.

In Section \ref{Norms Determined by Partitions and Weights},
we present well known examples of complemented subspaces 
of $L_p$ with norm given by partitions and weights.
We discuss some natural conditions on the families and in particular
normalization of the basis by inclusion of discrete 
partitions.
We also prove that
the natural sums of such Banach spaces are stable under these norms, i.e.,
have norms which are also given by partitions and weights.

In Section \ref{Embedding into L_p},
we first observe that if the norm on a space $X$ is given by finitely 
many partitions and weights, then $X$ is isomorphic to a subspace of 
$L_p$. Then we give the definition of an envelope
norm  and we prove the
existence of the envelope norm generated by a family of partitions
and weights. We also give a lower bound on a norm which is necessary
for a space to be isomorphic to a subspace of $L_p$. Finally we show
that if a space with norm given by partitions and weights is
isomorphic to a subspace of $L_p$, then its norm is equivalent to
the associated envelope norm. 

In Section \ref{Distance}, we construct examples which demonstrate the
difference between a norm given by partitions and weights and the
corresponding envelope norm. As a consequence we obtain an estimate
of the distance
between a certain Banach space $Y_n$ with norm given by partitions and weights
and $\otimes_{k=1}^n{X_p}$. 
Finally we give an example of a Banach space with norm given by
partitions and weights which is 
not isomorphic to a subspace of $L_p$ by applying the results from Section
3. Thus we get that not every sequence space with norm given by
partitions and weights is a ${\mathcal{L}}_p$ space. 

In the last sections we pose some questions for further study. In
particular we discuss the Bourgain, Rosenthal and Schechtman
construction and define spaces $X_{p}^{\alpha}$ with norm given by 
partitions and weights which are natural candidates for sequence space
realizations of the spaces $R_{p}^{\alpha}$.

We will use standard terminology and results in Banach theory as may 
be found in the books {\bf{[L-T-1]}},{\bf{[L-T-2]}} and {\bf{[J-L]}}.
Many results on subspaces of $L_p$ may be found in the exposition 
{\bf{[A-O]}} and its references.



\section{Norms Determined by Partitions and Weights}
\label{Norms Determined by Partitions and Weights}
In this section, we examine some examples of complemented subspaces 
of $L_p$ in order to motivate the idea of a norm given by partitions 
and weights. Then we develop the formal definition of a norm given 
by an admissible family of partitions and weights. Finally we give 
some results about sums of spaces with these norms.\\

In the following we will see that many well known complemented
subspaces of $L_p$ have equivalent norms of the form defined in
the Introduction. Here it is sometimes convenient to take partitions and 
weights defined on sets other than $\mathbb{N}$.
For each example we will have a family of partitions $(P_k)$ of 
${\mathbb{N}}^m$ for some $m$ and weights
$(W_k)$ for $k$ in some index set $K$.
\begin{example} Examples with one partition and weight.
 $K=\{1\}$.
    \begin{enumerate}
    \item If $P=\left\{ \{i\}: i\in \mathbb{N} \right\}$
and $W=(w_n)$ is any sequence of positive numbers, then 
$X\thicksim \ell_p$ since
$$
\|(x_n)\|=\|(x_n)\|_{P,W}
=\left( \sum_{n=1}^{\infty} (|x_n|^2w_{n}^{2})^{\frac p
  2}\right)^{\frac 1 p}
=\left( \sum_{n=1}^{\infty} |x_n|^pw_{n}^{p}\right)^{\frac 1 p}.
$$
    \item If $P=\{\mathbb{N}\}$ and $W=(w_n)$ is any sequence of 
positive numbers, then  $X\thicksim \ell_2$ since
$$
\|(x_n)\|=\left( (\sum_{n=1}^{\infty} |x_n|^2w_{n}^{2})^{\frac p
  2}\right)^{\frac 1 p}
=\left( \sum_{n=1}^{\infty} |x_n|^2w_{n}^{2}\right)^{\frac 1 2}.
$$

    \item If the index set is $\mathbb{N} \times \mathbb{N}$, 
the partition  $P=\left\{ \{n\} \times {\mathbb{N}}: n\in {\mathbb{N}}
\right\}$, and $W=(w_{n,m})_{n,m \in \mathbb{N}} $,
 then $X\thicksim \left(\sum \ell_2\right)_{\ell_p}$ since\\
$$
\|(x_n)\|=\left( \sum_{n=1}^{\infty} 
\left(\sum_{m=1}^{\infty}|x_{n,m}|^2w_{n,m}^{2}\right)^{\frac p
  2}\right)^{\frac 1 p}.
$$
    \end{enumerate}
\end{example}
\vskip .3 cm
\begin{example}\label{R-space}
 Examples with two partitions and weights.
 $K=\{1,2\}$
  \begin{enumerate}
    \item If $P_1=\big\{ \{n\} \big\}$ with weight $W_1=(1)$ and
 $P_2=\{{\mathbb{N}}\}$ with weight $W_2=(w_n)$, then $X$ is the
 space $X_{p, W_2}$, defined by Rosenthal, with norm
 $$\| (a_i)\|=\max \left\{ \left(\sum |a_n|^p\right)^{\frac{1}{p}}, 
\left(\sum |w_na_n|^2\right)^{\frac{1}{2}}\right\}.$$
  Rosenthal, {\bf{[R]}}, 
proved the following:
          \begin{enumerate}
          \item If $\begin{displaystyle}\inf_n w_n > 0
 \end{displaystyle}$,
 then $X_{p,W_2}\thicksim \ell_2$.
          \item If $\sum w_{n}^{\frac{2p}{p-2}} < \infty$,
 then $X_{p,W_2}\thicksim \ell_p$.
          \item If there is some $\epsilon >0 $ for which
 $\{n: w_n \ge \epsilon\}$ and
  $\{n: w_n < \epsilon\}$ are both infinite and for which
 $\begin{displaystyle}\sum_{w_n <\epsilon} w_{n}^{\frac{2p}{p-2}} < 
\infty\end{displaystyle}$,
 then $X_{p,W_2} \thicksim \ell_2 \oplus \ell_p$.
          \item $\text{For}\;\text{ each}\;\epsilon > 0, 
\begin{displaystyle}\sum_{w_n< \epsilon}
  w_{n}^{\frac{2p}{p-2}}=\infty.\end{displaystyle}\phantom{xxxxxx}(*)$\\
If $W_2$ satisfies $(*)$, then $X_{p,W_2} \thicksim X_p$.
          \end{enumerate}

     \item  If $P_1=\left\{ \{(i,j)\} \right\}$ with weight $W_1=(1)$
 and
 $P_2=\left\{ \{n\}\times\mathbb{N}\right\}$ with weight
 $W_2=(w_{n,m})$ where
 $w_{n,m}=\left(\frac{1}{n}\right)$ 
for all $n$,$m$, 
 then $X\thicksim \left(\sum_n X_{p, (\frac{1}{n})}\right)_p$.\\ 
 $X_n = X_{p,(\frac{1}{n})}$ is isomorphic to $\ell_2$ and
 $\sup_{n\in \mathbb{N}} d(X_n, \ell_2) = \infty$, so 
 $X \thicksim B_p$, as defined by Rosenthal.
    \end{enumerate}
\end{example}

\begin{example} An example with four partitions and weights.
$K=\{0,1,2,3\}$.\\
Let $i$ represent the first index and $j$
represent the second index in the set $\mathbb{N}\times \mathbb{N}$.
Assume the sequences $(w_i)$ and $(w'_j)$ satisfy $(*)$. 
\begin{alignat*}{3}
\text{Let } P_0 &=\{\mathbb N \times \mathbb N\}  &&\qquad \text{ with weight}
&\qquad W_0 &=(w_i\cdot w'_j)\\
P_1&=\left\{ \{n\} \times \mathbb{N}:n\in\mathbb N\right\}
&&\qquad \text{ with weight} & W_1&=(1\cdot w'_j)\\
P_2&=\left\{{\mathbb{N}}\times \{n\} :n\in\mathbb N\right\}
 &&\qquad \text{ with weight} & W_2&=(w_i \cdot 1)\\
P_3&= \{\{(i,j)\}:i,j\in \mathbb N\} &&\qquad \text{ with weight} &
 W_3&=(1\cdot 1)
\end{alignat*}
Then this is Schechtman's basic example, {\bf{[S]}},
$X \thicksim X_p \otimes X_p$, with norm 
\begin{multline}
\max \left\{ \left(\sum_{i,j} |a_{i,j}|^2w_{i}^{2}
{w'}_{j}^{2}\right)^{\frac{1}{2}}, \left(\sum_i\left(\sum_j |a_{i,j}|^2
{w'}_{j}^{2}\right)^{\frac{p}{2}}\right)^{\frac{1}{p}},\right .\\
\left .  \left( \sum_j \left( \sum_i
|a_{i,j}|^2 w_{i}^{2} \right)^{\frac{p}{2}} \right)^{\frac{1}{p}},
\left(\sum_{i,j} |a_{i,j}|^p \right)^{\frac{1}{p}} \right\}
\thickapprox\left\|\sum_{i,j} a_{i,j}
(x_i \otimes y_j)\right\|_{L_p(I\times I)}
\end{multline}
\end{example}

This example can be generalized by using the index set
 ${\mathbb{N}}^n$. If $|K|=2^n$ and the partitions and weights are chosen
 in a manner similar to the above, then $X\thicksim \otimes^n_{k=1}X_{p}$.
Thus we get isomorphs of Schechtman's examples. Additional detail
about these examples is contained in Section \ref{Distance}.

We will now give a general definition of a space with norm given by
partitions and weights.  Below $A$ is any countable set. 

\begin{definition}\label{norm}
Let $P=\{ N_i\}$ be a partition of $A$ and 
 $ W : A \rightarrow { ( 0, 1 ]}$  be a function, the {\it weights}. 
Let $x_j \in \mathbb{R}$ for all $j\in A$.
Define \\
 $$\begin{displaystyle} \|(x_j)_{j\in A}\|_{P,W} = \left(\sum_{i} \left (
 \sum_{j\in {N_i}} {x_j}^2 w_{j}^{2}\right
 )^{\frac{p}{2}}\right)^{\frac{1}{p}}
 \end{displaystyle}$$\\
Suppose that $(P_k, W_k)_{k\in K}$ is a family of pairs of partitions
 and functions as above. Define a (possibly infinite)
norm on the real valued function on 
$A$, $(x_i)_{i \in A},$ by
 $ \|(x_i)\| = \begin{displaystyle} \sup_{k\in K} \|(x_i)\|_{(P_k,
 W_k)} \end{displaystyle} $ 
and let $X$ be the subspace of elements of finite norm. In this case
we say that $X$ has a {\it norm given by partitions and weights.}
\end{definition}
\begin{remark}
Because of the nature of this norm, $X$ will have a natural
 unconditional basis. Thus this approach to describing the
 complemented subspaces of $L_p$ is limited to complemented
 subspaces of $L_p$ with unconditional basis. At this time,
 no complemented subspace of $L_p$ without unconditional basis is known.
\end{remark}

\begin{proposition}
Suppose $X$ has a norm given by partitions and weights. Then $X$ is a Banach
space.
\end{proposition}
We leave the straight-forward proof to the reader.

\begin{proposition}\label{1-partition} Suppose $X$ is a 
 Banach space with norm
 given by one partition and weight. Then  $X \thicksim \ell_p$,
 $X \thicksim \ell_2$,
$ X \thicksim \ell_2 \oplus \ell_p$, or  $X \thicksim \begin{displaystyle}
 \left ( {{\sum}^{\oplus} \ell_2}\right )_{\ell_p}  \end{displaystyle}$.
\end{proposition}

Notice that these are the spaces given in Example 2.1 and 
their direct sums. 
The proof is a routine computation after normalization of the basis.

Since normalization of the basis is an important first step
to understanding the spaces, we now 
introduce admissible families of partitions and weights to incorporate
this and some other properties.

\begin{definition}
The partition of $A$, $\{ \{a\}: a\in A\}$, will be called the {\it
discrete}
 partition. The partition of $A$, $\{ A\}$,
 will be called the {\it indiscrete} partition. 
\end{definition}

\begin{definition}
A family of partitions and weights is called {\it admissible} if it contains
the discrete partition with the trivial weight $(w(a))_{a\in A}=(1)$ and
the indiscrete partition with some weight.
\end{definition}

The discrete partition is included to force the natural coordinate
basis to be normalized. This requirement is not really a restriction
because every normalized unconditional basic sequence in $L_p$
has a lower $\ell_p$ estimate. (See the Preliminaries section of {\bf [A-O]}.)
The indiscrete partition gives a candidate
for a natural $\ell_2$ structure on the space $X$. Because we are
concerned with embedding these spaces into $L_p, \; p>2$, 
there always must be some $\ell_2$ structure on the space.

Notice that in the previous examples, Rosenthal's space and
Schechtman's space have norms given by admissible families 
of partitions and weights. Each of the other cases can be 
equivalently renormed using an admissible family of partitions 
and weights. Unless otherwise noted we will assume from now on
that a Banach space $X$ with norm given by partitions and weights
is actually given by an admissible family of partitions and weights.

Next we are going to prove some stability results
for sums of spaces when the spaces are equipped with these norms.

\begin{definition}\label{(p,2)-sum}
Let $(X_n)$ be a sequence of subspaces of 
$L_p(\Omega, \mu)$ for some 
probability measure $\mu$, and let $(w_n)$ be a sequence of real
numbers, $0< w_n \le 1$. For any sequence $(x_n)$ such that $x_n \in
X_n$ for all $n$, let
$$\|(x_n)\|_{p,2,(w_n)}=\max \left\{\left(\sum
\|x_n\|_{p}^{p}\right)^{\frac 1 p}, \left(\sum \|x_n\|_{2}^{2} 
w_{n}^{2}\right)^{\frac 1 2}\right \}$$
and let
\begin{align*}
X&=\left(\sum X_n\right)_{(p,2,(w_n))}
=\{(x_n): x_n \in X_n
\;\text{for}\;\text{all}\;\|(x_n)\|_{p,2,(w_n)}<\infty\}
\end{align*}
We will say that $X$ is the {\it $(p,2,(w_n))$-sum of $\{X_n\}$.}
\end{definition} 

Let $A$ be a countable set and let $(X_a)_{a\in A}$ be a 
family of Banach spaces of functions defined on sets $(B_a)_{a\in A}$,
respectively. That is, for each $a\in A$, $X_a$ has a norm given by
a family of partitions of $B_a$ and weights on $B_a$. Let $I_a$ denote the
index set of the corresponding family for $X_a$. For each $i(a)\in
I_a$, let $P^{a, i(a)}$ be a partition of $B_a$
and $W^{a, i(a)}$ be a weight function, i.e.,
$W^{a, i(a)}: B_a \rightarrow (0, 1]$. For each $a\in A$ and $i(a)\in
I_a$, the norm on $X_a$
with respect to $P^{a, i(a)}, W^{a, i(a)}$ is given by
$$\|(x_{a,b})_{b\in B_a}\|_{P^{a,i(a)},
 W^{a,i(a)}} =\left ( \sum_{Q \in P^{a, i(a)}} \left(\sum_{b \in Q }
 (x_{a,b})^2(w^{a,i(a)}(b))^2\right)^{\frac{p}{2}}\right )^{\frac{1}{p}}
$$
For each $a$, we assume that there is one distinguished indiscrete partition
and weight. We will denote the index of this partition and weight as $(\;)$.
Let $P^{a,(\;)} = \{ B_a \}$, and $W^{a,(\;)}$ be the associated weight.
For each $a$, define
$\begin{displaystyle} \label{2 norm} 
\|(x_{a,b})_{b\in B_a}\|_2 
= \left(\sum_{b \in B_a } (x_{a,b})^2(w^{a,(\;)}(b))^2\right)^{\frac{1}{2}} 
\end{displaystyle}$.
Suppose that for the index set $A$, we have an associated weight function
$W: A\rightarrow (0, 1]$.
Let $(\sum_{a\in A} X_a)_{p,2,W} $
be defined on $B= {\coprod}_{a\in A} B_a$ as above using the norm
$\|(x_{a,b})_{b\in B_a}\|_2 $
as the $\|\cdot\|_2$ in the definition.
Let $I = {\prod}_{a\in A} I_a \cup \{(\phantom{x})\}$. Let
$(i(a))_{a\in A} \in I$.
Then there is a natural partition of $B$ and weight on $B$ given by
$P_{(i(a))}=\left\{\{a\}\times P : P \in P^{a, i(a)}, a\in A\right \}$  and 
$W_{(i(a))}=(w_{b}^{a,i(a)})_{b\in B_a, a\in A}$.
We define as a special case the partition and weight for $(\phantom{x}
 )$ as
 $P_{(\phantom{x} )}=\{  {\coprod}_{a\in A} B_a\}$ and
 $W_{(\phantom{x} )}=(W(a)w^{a, (\phantom{x} )}(b))_{b\in B_a, a \in A}$.\\
If we expand the definition of the norm we have
\begin{eqnarray}
&&\|(x_{a,b})_{a\in A, b\in B_a}\|_{p,2,W}\nonumber\\
&=&\max \left \{ \left(\sum_{a\in A} \|(x_{a,b})_{b\in
  B_a}\|_{X_a}^{p}\right)^{\frac{1}{p}},
 \left(\sum_{a\in A} \|(x_{a,b})_{b\in B_a}\|_{2}^{2}
(W(a))^2\right)^{\frac{1}{2}}\right \}\nonumber\\
&=&\max \left \{ \left(\sum_{a\in A} \sup_{i(a)\in I_a}
 \{\|(x_{a,b})_{b\in B_a}\|_{P^{a,i(a)}, 
W^{a, i(a)}}^{p}\}\right)^{\frac{1}{p}},\right .\nonumber\\
&&\phantom{xxxxxxxxxxxxxxxxxxxx}\left .\left(\sum_{a\in A} 
W(a)^2\sum_{b\in B_a}
(w^{a, ( )}(b))^2 |x_{a,b}|^2\right)^{\frac{1}{2}}\right \}\label{p,2,W-sum}
\end{eqnarray}

Notice that for each $a\in A$, we take a supremum over $I_a$,
then we take summation of those supremums, and finally we take the
maximum of two sums. If we consider the index $(i(a))$ which for each
$a$ approximates the supremum, it is one element in $I$. So instead
of taking the maximum over each $I_a$, we can compute the norm for
each index in
$I$, and then take supremum of them only once. Hence the norm becomes
$$ \left\|(x_{a,b})_{a\in A, b\in B_a}\right\|_{p,2,W}=\sup_{(i(a))\in I}
 \left\|(x_{a,b})_{ b\in B_a}\right\|_{P_{(i(a))}, W_{(i(a))}}$$
This gives us the following result:
\begin{proposition}\label{p,2 sum}
Let $(X_a)_{a\in A}$ be a family of Banach spaces each with norm given
 by partitions and weights. Then the norm of the space $(\sum_a
 X_a)_{p,2, W}$ can also be expressed as a norm given by partitions
 and weights. In other words, the class of spaces with norm given by
 partitions and weights is stable under $p,2$ sums.
\end{proposition}

\begin{corollary}
Let $(X_a)_{a\in A}$ be a family of Banach spaces with norm given
 by an admissible family of partitions and weights. Then the norm 
of the space $\left(\sum_a X_a\right)_{\ell_p}$ can also be expressed with
 partitions and weights.
\end{corollary}

Proof: Choose $W$ such that $\begin{displaystyle}\sum_{a\in A} 
W(a)^{\frac{2p}{p-2}} < 1 \end{displaystyle}$ in the previous result.


\section{Embedding into $L_p$}
\label{Embedding into L_p}

In this section, we first show that any space $X$ with norm
given by finitely many partitions and weights is
isomorphic to a subspace of $L_p$.
Then we give the definition of an envelope
norm. We prove the
existence of the envelope norm generated by a family of partitions
and weights. We also give a lower bound on a norm which is necessary
for a space to be isomorphic to a subspace of $L_p$. Finally we show
that if a space with norm given by partitions and weights is
isomorphic to a subspace of $L_p$, then its norm is equivalent to
the natural envelope norm.

\begin{proposition}
Any sequence space $X$ with norm given by finitely many partitions and
weights is isomorphic to a subspace of $L_p$.
\end{proposition}

Proof: Let $X$ be the sequence space with partitions and weights
$(P_n, W_n)_{n=1}^{N}$. Let $X_n$ be the space of sequences with norm
given by one partition
and weight $(P_n, W_n)$, $1\leq n \leq N$. Then it is easy to see that, 
$$\begin{displaystyle}\left(\sum_{n=1}^{N} X_n\right)_{\ell_{\infty}}\thicksim
\left(\sum_{n=1}^{N} X_n\right)_{\ell_p}\end{displaystyle}$$
Take the isometric embedding from $X$ into 
$\begin{displaystyle}\left(\sum_{n=1}^{N}
X_n\right)_{\ell_{\infty}}\end{displaystyle}$ 
defined by $x\mapsto (x)_{n=1}^{N}$. Since 
 for each $n$, $n=1,2,...,N$, $X_n$ is
isomorphic to a complemented subspace of $L_p$
by Proposition~\ref{1-partition},
then  $\begin{displaystyle}\left(\sum_{n=1}^{N}
X_n\right)_{\ell_p}\end{displaystyle}$ is isomorphic to a complemented
subspace of $L_p$. Hence $X$ is isomorphic to a subspace of $L_p$.
\qed \\

This trivial approach fails if there are infinitely many partitions
and weights. We do not know any general sufficient condition on the
partitions and weights to guarantee that the space is isomorphic to a 
subspace of $L_p$. By using the fact that these spaces have
unconditional basis, we can give a useful necessary condition.\\

\begin{definition}\label{Def-enve-norm} Let $X= \{ (a_b)_{b\in B}\}$ 
be a Banach space defined 
on a countable set $B$ with norm given by a set of partitions and
 weights
 ${\mathcal{P}}=\{(P^i, W^i): i\in K\}$. 
$\mathcal{P}$ satisfies 
the {\it envelope property}
if and only if for any
 partition $Q$ of $B$, and any function $\phi: Q \rightarrow
 {\mathcal{P}}$, $\phi(q)= (P^{i(q)}, W^{i(q)})$ for all $q\in Q$, 
the partition and weight
 $(P_0, W_0)$ belongs to $ \mathcal{P}$ where 
$P_0=\{ K : K = q \cap K_{i(q)}\neq \varnothing \phantom{x}\text{for}
\phantom{x} \text{some}\phantom{x} q\in Q, \text{some}
\phantom{x} K_{i(q)}\in P^{i(q)}\}$ and 
$W_0= (w_{b}^{i(q)})_{b\in q, q\in Q} \phantom{x}\text{where}
\phantom{x} W^{i(q)}=(w_{b}^{i(q)})_{b\in B}$.\\
In this case we will say that
$\|| \cdot \|| = \sup_{i\in K}
 \| \cdot \|_{P^i,W^i}$ is an {\it envelope norm.}
\end{definition}
\begin{example} \label{Xp} Let $X_p$ be the Rosenthal's space with norm 
$$\| (a_i)\|=\max \left\{ \left(\sum |a_n|^p\right)^{\frac{1}{p}},
 \left(\sum |w_na_n|^2\right)^{\frac{1}{2}}\right\} $$
where $(w_n)$ satisfies (*).
Let $P_1=\left\{ \{n\} \right\}$ with weight $W_1=(1)$ and
 $P_2=\{{\mathbb{N}}\}$ with weight $W_2=(w_n)$.
 Then ${\mathcal{P}}=\{ (P_1, W_1), (P_2, W_2)\}$ defines the norm on
 $X_p$. It is easy to see that $\mathcal{P}$ does not have the
 envelope property. To get a family of partitions and weights which
 has the envelope property we need to add all of the possible
 combinations of the given two.\\
Let $\mathcal{Q}$ be the set of all partitions on $\mathbb{N}$.
 Let $Q\in \mathcal{Q}$ and $T: Q \rightarrow \mathcal{P}$. Define
\begin{multline*} P(Q, T)=\{ K: K=\{ n\} \text{ if $n\in q$ and } 
T(q)=(P_1, W_1) \text{ for some $q\in Q$}\}\\
\cup\{K: K=q \text{ if }
 T(q)=(P_2, W_2)\text{ for some $q\in Q$}\}
\end{multline*}
and
\begin{equation*}W(Q,T)= (w(n))_{n\in q, q\subset \mathbb{N}}\text{ where
}\\w(n)=\begin{cases} 
1& \text{if $n\in q$, $T(q)=(P_1, W_1)$ }\\
w_n& \text{if $n\in q$, $T(q)=(P_2, W_2)$}.
\end{cases}
\end{equation*}
Then an equivalent envelope norm is defined by 
$\sup_{(P,W)\in \tilde{\mathcal{P}}}\|\cdot\|_{(P,W)}$ where\\
$\tilde{\mathcal{P}}=\{(P(Q,T),W(Q,T)): Q \in{\mathcal{Q}},
 T: Q\rightarrow \mathcal{P}\}$.
\end{example}

\begin{remark} In the case of the norm  of $X_p$ in Example \ref{Xp}
the envelope norm can be written in the form
$$\max_{q\subset \mathbb N} \left( \sum_{n\in q} |a_n|^p +(\sum_{n \notin q}
|a_n|^2 w_n^2)^{p/2} \right )^{1/p}.$$
\end{remark}

Next we show that by generalizing the construction above that there is 
a natural envelope norm associated to each norm given by partitions
and weights.

\begin{proposition}
\label{existence}
Suppose $X$ is a Banach space defined on a countable set $B$ with norm
 given by a family $\mathcal{P}$ of partitions and weights. Then there
 exists a natural family of partitions and weights
 $\tilde{\mathcal{P}}$, (defined below), such that
 $\||\cdot\||= \sup_{(P,W)\in \tilde{\mathcal{P}}}
 \|\cdot\|_{P,W}$
 is an envelope norm.
\end{proposition}

Proof: Let $\mathcal{Q}$ be the set of all partitions of $B$.
Let ${\mathcal{P}}=\{(P_i, W_i): i\in K\}$ be the given family of 
partitions and weights for $X$.
Let $Q\in \mathcal{Q}$. Let $T$ be a map from $Q$ into 
 $\mathcal{P}$ denote  $T(q)$ by $(P^{i(q)}, W^{i(q)})$ for all $q \in Q$.
 Define 
$P(Q,T)=\{K: K=q\cap p\neq \varnothing, q\in Q, p \in P^{i(q)}\}$ and
$W(Q,T)=(w^{i(q)}(b))_{b\in q, q\in Q} \; \text{where}\; 
 W^{i(q)}=(w^{i(q)}(b))_{b\in B})$.
Let $\tilde{\mathcal{P}}=\{ (P(Q,T), W(Q,T)): Q\in {\mathcal{Q}},\;
 T: Q\rightarrow \mathcal{P}\}$.
Define a norm on $X$ as $\|| (x_i)\|| = \sup_{(P,W)\in 
\tilde{\mathcal{P}}} \|(x_i)\|_{P,W}$.
We claim that $\tilde{\mathcal{P}}$ has the envelope property and thus
$\||\cdot\||$ is an envelope norm.

Let $\bar{Q}$ be any partition of $B$.  Let $S$ be any map from
 $\bar{Q}$ into $\tilde{\mathcal{P}}$, i.e.,
$$S(\bar{q})=\left(P(Q_{\bar{q}}, T_{\bar{q}}),W(Q_{\bar{q}},
 T_{\bar{q}})\right)$$
 for all  $\bar{q} \in \bar{Q}$. 
For any $\bar{q}\in \bar{Q}$, let
 $T_{\bar{q}}(q_0)=(P^{i(\bar{q}, q_0)}, W^{i(\bar{q}, q_0)})$ 
for all $q_0\in Q_{\bar{q}}$ and 
let $\overline{\overline{Q}}=\{\overline{\overline{q}}\neq \varnothing :
 \overline{\overline{q}}=q_0 \cap \bar{q}, \bar{q}\in
 \bar{Q}, q_0\in Q_{\bar{q}}\}$.
Because $\bar{Q}$ and $Q_{\bar{q}}$ are partitions,
 $\overline{\overline{q}}$ uniquely determines $\bar{q}
 \in \bar{Q}$ and $q_0 \in Q_{\bar{q}}$ such that 
$\overline{\overline{q}}=q_0 \cap \bar{q}$.
From Definition \ref{Def-enve-norm} we have
$P_0=\{\bar{K}\neq \varnothing: \bar{K}=\bar{q}\cap 
\bar{K}_{\bar{q}}, \bar{q}\in \bar{Q}, \bar{K}_{\bar{q}}
\in P(Q_{\bar{q}},T_{\bar{q}})\}$,
which is exactly what the definition above gives 
 for the partition of $B$ determined by 
$\bar{Q}$ and $S$, $P(\bar{Q}, S)$. Thus
\begin{eqnarray}
P_0&=&\bar{P}(\bar{Q}, S)\nonumber\\
&=&\{\bar{K}: \bar{K}=\bar{q}\cap 
\bar{K}_{\bar{q}}\neq \emptyset, \bar{q}\in \bar{Q}, \bar{K}_{\bar{q}}
\in P(Q_{\bar{q}},T_{\bar{q}})\}\nonumber\\
&=&\{\bar{K}: \bar{K}=\bar{q}\cap(q_0\cap p)\neq \emptyset, 
\bar{q}\in \bar{Q},
 q_0 \in Q_{\bar{q}}, p\in P^{i(\bar{q},q_0)}\}\label{3.1}\\
&=&\{\bar{K}: \bar{K}=(\bar{q}\cap q_0)\cap p\neq \emptyset, \bar{q}\cap q_0\in
 \overline{\overline{Q}}, p\in P^{i(\bar{q}, q_0)}\}\label{3.2}\\
&=&\{\bar{K}: \bar{K}=\overline{\overline{q}}\cap p\neq \emptyset,
 \overline{\overline{q}}=\bar{q}\cap q_0\in \overline{\overline{Q}},
 p\in P^{i(\bar{q}, q_0)}\}\label{3.3}
\end{eqnarray}
where (\ref{3.1}) follows from the definition of $P(Q_{\bar{q}},
T_{\bar{q}})$, (\ref{3.2}) by the definition of $\overline{\overline{Q}}$, and
(\ref{3.3}) by the uniqueness of $\bar{q}$ and $q_0$.
Define $\overline{\overline{T}}: \overline{\overline{Q}}
 \rightarrow \mathcal{P}$ by 
$ \overline{\overline{T}}(\overline{\overline{q}})
=(P^{i(\bar{q}, q_0)},W^{i(\bar{q}, q_0)})$
where $\overline{\overline{q}}=\bar{q}\cap q_0$, $q_0 
\in Q_{\bar{q}}$, $\bar{q}\in \bar{Q}$. Then we have shown that
 $\bar{P}(\bar{Q},
 S)=P(\overline{\overline{Q}},\overline{\overline{T}})$.\\

Because $T_{\bar{q}}(q)=\left( P^{i(\bar{q},q)},
  W^{i(\bar{q},q)}\right)=\left( P^{i(\bar{q},q)},
  (w^{i(\bar{q},q)}(b))_{b\in B}\right)$, \\
$$S(\bar{q})=\left( P(Q_{\bar{q}}, T_{\bar{q}}),W(Q_{\bar{q}}, T_{\bar{q}})
  \right)=\left( P(Q_{\bar{q}},
    T_{\bar{q}}),(w^{i(\bar{q},q)}(b))_{b\in q,q\in
      Q_{\bar{q}}}\right).$$\\
Suppose $W_0=(w_b)_{b\in B}$. If $\bar{q}\in \bar{Q}$ and $b\in
\bar{q}$, then as in Definition \ref{Def-enve-norm},
$w_b=w^{i(\bar{q},q_0)}(b)$
where $b\in q_0$ and $q_0 \in Q_{\bar{q}}$.
Hence for $b\in \overline{\overline{q}}=\bar{q}\cap q_0$, 
$w^{i(\bar{q},q_0)}(b)$ is also the choice specified by 
$\overline{\overline{T}}(\overline{\overline{q}})$. 
Hence $W_0=W(\overline{\overline{Q}}, \overline{\overline{T}})$.
So $(P_0, W_0)\in \tilde{\mathcal{P}}$. \qed

\begin{corollary} If $X$ has a norm defined by a finite number of
  partitions and weights, then there is an equivalent envelope norm on
$X$.
\end{corollary}

The next result follows by simply checking that all of the partitions
and weights in the construction are required for the envelope property.\\

\begin{proposition}
\label{minimal}
Suppose $\mathcal{P}$ is a family of partitions and weights on $B$.
 Then $\tilde{\mathcal{P}}$ as in Proposition $\ref{existence}$ 
is the minimal family of partitions and weights on $B$ containing
 $\mathcal{P}$ and satisfying the envelope property.
\end{proposition}

\begin{corollary}
Let $\mathcal{P}$ be a non empty family of partitions and weights on
 $B$.
 Let $(\mathcal{P_{\lambda}})_{\lambda \in \Lambda}$ be any chain of
 families of partitions and weights on $B$ such that each
 $\mathcal{P_{\lambda}}$ contains $\mathcal{P}$ and satisfies envelope
 property. Then $\cap_{\lambda \in \Lambda} \mathcal{P}_{\lambda}$
 satisfies the envelope property.
\end{corollary}

Our purpose in introducing envelope norms is to show that the 
envelope norm is related to a property of subspaces of $L_p$ 
with unconditional basis that we will now explain.

Let $X$ be a Banach space defined on $B$ with a norm given 
by partitions and weights. Let $\phi$ be a one-to-one map 
from $\mathbb{N}$ onto $B$ such that  $x_n=e_{\phi(n)}$
where $(e_b)_{b\in B}$ is the natural unit vector basis of $X$.
Then $(x_n)$ is an unconditional basis for $X$. Let
$\begin{displaystyle}x=\sum_{n=1}^{\infty} a_n x_n\end{displaystyle}$ 
for some $(a_n)$.
Let $Q$ be any partition of $B$.  Let $\{F_k\}_{k=1}^{\infty}$
be the corresponding partition of $\mathbb{N}$, i.e., $\phi(F_k)=q$
for some $q\in Q$. Then
$$ x=\sum_{k=1}^{\infty}\sum_{n\in F_k} a_n x_n
 = \sum_{k=1}^{\infty} z_k 
= \sum_{q\in Q} z_{q}'
$$
where 
$$
z_k = \sum_{n\in F_k} a_nx_n=\sum_{n\in F_k} a_ne_{\phi(n)}
=\sum_{b\in q=\phi(F_k)} a_{\phi^{-1}(b)}e_b=z_{\phi(F_k)}'=z_{q}'
$$
Since $(x_n)$ is an unconditional basis and $(z_k)$ (hence $(z_{q}')$) 
is a block of $(x_n)$, $(z_k)$ is an unconditional basic sequence
with unconditional constant 1.
\begin{remark}
We are abusing the terminology a little here. It would be more correct to
use partitions $\{F_k\}$
containing only finite subsets of $\mathbb N$ and the
map $\phi$ which satisfies the property that for all $k$,
if $n \in \phi^{-1}(F_k)$ and $m\in
\phi^{-1}(F_{k+1})$ then $n<m$.
\end{remark}

In the lemma below we use the notation introduced in Proposition 
\ref{existence}.

\begin{lemma}\label{lower estimate}
Let ${\mathcal{P}}=\{(P^i, W^i): i\in K\}$ be a family of 
partitions and weights on $B$. 
Let $X$ be the corresponding Banach space defined on $B$. 
 If $X$ is isomorphic to a subspace of $L_p$, then there exists 
a constant $C$, depending only on the Banach-Mazur distance of $X$ to a
 subspace of $L_p$, such that for any partition $Q$ of $B$ and
 any map $T: Q\rightarrow \mathcal{P}$, $\|x\|\ge C\|x\|_{(P(Q,T),W(Q,T))}$
 where $T(q)=(P^{i(q)}, W^{i(q)})$.
\end{lemma}

Proof: Let $\phi: {\mathbb{N}} \rightarrow B $ as above and 
$T: Q \rightarrow \mathcal{P}$ such that $T(q)=(P^{i(q)}, W^{i(q)})$.
 If $X$ is isomorphic to a subspace of $L_p$, with isomorphism
 $R$, then $(Rz_k)$ (hence $(Rz_{q}')$) is block of $(Rx_n)$
 which is an unconditional basic sequence in $L_p$ with constant
 $\lambda$. So 

\begin{align}
\|x\|_X&=\|\sum_{k=1}^{\infty} z_k\|
\ge \|R\|^{-1}\|\sum_k Rz_k\|_{L_p}\nonumber\\
&\ge \|R\|^{-1}\lambda^{-1}\left(\sum_k\|Rz_k\|_{p}^{p}\right)^
{\frac{1}{p}}\label{a-2}\\ 
&\ge \|R\|^{-1}\lambda^{-1}\left(\sum_k \frac {\|z_k\|_{X}^{p}}
{\|R^{-1}\|^p}\right)^{\frac{1}{p}}
=\frac{\lambda^{-1}}{\|R\|\|R^{-1}\|}\left(\sum_{q\in Q}
\|z_q\|_{X}^{p}\right)^{\frac{1}{p}}\nonumber\\
&\ge C\left(\sum_{q\in Q}\sum_{r\in P^{i(q)}} 
\left(\sum_{\phi(n) \in {r\cap q}}
 |a_{n}|^2(w_{\phi(n)}^{i(q)})^2\right)^{\frac{p}{2}}\right)^
{\frac{1}{p}}\label{a-4}\\
&=C\left(\sum_{\bar{q}\in P(Q,T)}\left(\sum_{\phi(n)\in \bar{q}} 
 |a_{n}|^2(w_{\phi(n)}^{i(q)})^2\right)^{\frac{p}{2}}\right)^
{\frac{1}{p}}
=C\|x\|_{P(Q,T), W(Q,T)}\nonumber
\end{align}
where
(\ref{a-2}) follows from the standard lower $\ell_p$ estimate,
 {\bf{[A-O]}} and  (\ref{a-4}) is true since 
$ z_q=\sum_{n\in F_k, \phi(F_k)=q} a_n x_n$, and 
$\|z'_q\|_X \ge \|z'_q\|_{P^{i(q)},W^{i(q)}}$. 
 In (\ref{a-4}), $q$ is the unique element of $Q$ such that
 $\bar{q}\subset q$.
\qed

\begin{theorem}\label{equivalent norm}
Suppose $X$ has a norm given by a family $\mathcal{P}$ of partitions and
 weights and $X$ is isomorphic to a subspace of $L_p$. Then
 there is an envelope norm $\||\cdot\||$ such that
 $\||\cdot \|| \sim \|\cdot\|_X$
\end{theorem}

Proof: If we take a supremum over all the choices of $Q$
 and $T$ in Lemma \ref{lower estimate}, we have $\|x\|_X \ge C \||
 x\||$, where $\||\cdot\||$ is the envelope norm defined by 
$\tilde{\mathcal{P}}$ in Proposition \ref{existence}. 
On the other hand, 
since ${\mathcal{P}}\subset \tilde{\mathcal{P}}$,
 we get $\|x\|_X \le \|| x\||$. Hence  $\|x\|_X \sim \|| x\||$\qed

\begin{remark}: Because the natural basis of a space with norm given 
by partitions and weights is 1 unconditional, the unconditional
constant of the image of any block basis under an isomorphism $R$ is
at most $\|R\|\|R^{-1}\|$. Hence the constant $\lambda$ in the proof of Lemma 
\ref{lower estimate} and consequently the equivalence in Theorem 
\ref{equivalent norm} depend only on the distance to a subspace of $L_p$.
\end{remark}

\begin{proposition}\label{p,2 envelope}
Let $(X_a)_{a\in A}$ be a family of Banach spaces each with norm given
 by partitions and weights which satisfies the envelope property. Then the 
norm of the space $(\sum_a X_a)_{p,2, W}$ can also be expressed as a
 norm given by partitions and weights which also satisfies the
 envelope property.
\end{proposition}
\begin{proof}This follows from equation \ref{p,2,W-sum}. Indeed, if
$Q$ is any partition of $\coprod_{a\in A} B_a$ and \newline
$\phi:Q\rightarrow
\{(P_{(i(a))},W_{(i(a))}): (i(a)) \in I\}$ then for each $q \in Q$
let $\phi(q)=(P_{(i(a),q)},W_{(i(a),q)})$. Then
$(Q_a,W_a)$, where
$Q_a=\{q\cap q'\ne \empty: q \in Q, q'\in P_{i(a),q}\}$ and
$W_a=(w^{a,i(a,q)}_b)_{b\in B_a}$, must be one of the partitions and
weights in
$\{(P_{a,i(a)},W_{a,i(a)}):i(a)\in I_a\}$  since 
this family of partitions and weights has the envelope property.
\end{proof}

Next we consider the $L_p$ tensor product of spaces with norms given by
partitions and weights. Because the tensor product is only defined for
subspaces of $L_p$, we will assume that the two spaces are also isomorphic
to subspaces of $L_p$ and thus the defining families of partitions and
weights must have the envelope property.

\begin{proposition}\label{tensor}
Suppose that $\mathcal P$ and $\mathcal Q$ are families of partitions and
weights defined on sets $A$ and $B$, respectively, and having the envelope
property. Let $X$ and $Y$, be the corresponding spaces on $A$ and $B$ and
let $(x_a)_{a\in A}$ and $(y_b)_{b\in B}$, respectively, be the natural bases.
Suppose that $S$ is an isomorphism from $X$ into $L_p[0,1]$ and $T$ is an
isomorphism from $Y$ into $L_p[0,1]$. Then there is a constant $C>0$ such
that for any constants $(c_{a,b})_{a\in A,b \in B}$ with only finitely many
non-zero,
any $(P,(w_a))\in \mathcal P$ and any $(Q,(w_b'))\in \mathcal Q$,
$$ \|\sum c_{a,b} Sx_a \otimes Ty_b\|_p \ge C \left(\sum_{p\in P, q\in Q}
\left( \sum_{a\in p,b \in q} |c_{a,b}|^2 (w_a)^2 (w_b')^2\right )^{p/2}\right
)^{1/p}.$$
\end{proposition}

\begin{proof} With the given notation we can directly estimate the norm as
follows.
\[
\begin{split}
 \|\sum c_{a,b} Sx_a \otimes Ty_b\|^p_p 
&= \int_0^1 \int_0^1 |\sum
c_{a,b} Sx_a(s)Ty_b(t)|^p\, dt\, ds\\
&= \int_0^1 \int_0^1 \left|\sum_{b\in B} (\sum_{a\in A} c_{a,b}
Sx_a(s))Ty_b(t)\right|^p\, dt\, ds \\
&=\int_0^1 \left\|T(\sum_{b\in B} (\sum_{a\in A} c_{a,b} 
Sx_a(s))y_b\right\|^p_p \, ds\, \\
& \sim \int_0^1 \max_{(Q_1,(w_b''))\in \mathcal Q} \left(\sum_{q\in
Q_1}\left(\sum_{b\in q}
|\sum_{a\in A} c_{a,b} Sx_a(s)|^2 (w''_b)^2\right)^{p/2}\right)\, ds\\
\intertext{(since $T$ is an isomorphism)}\\
&\ge \int_0^1  \left(\sum_{q\in Q}\left(\sum_{b\in q}
|\sum_{a\in A} c_{a,b} Sx_a(s)|^2 (w'_b)^2\right)^{p/2}\right)\, ds\\
&\sim   \sum_{q\in Q} \int_0^1 \int_0^1 \left|\sum_{b\in
q}\sum_{a\in A} c_{a,b} Sx_a(s)w'_b r_b(u)\right|^p\,du \, ds \\
\intertext{(by Khinchin's inequality, where $(r_b)_{b\in B}$ is a sequence
of Rademacher functions.)}
&=   \sum_{q\in Q} \int_0^1 \left \|\sum_{a\in
A}(\sum_{b\in q}  c_{a,b}w'_b r_b(u))Sx_a(s)\right\|^p\,du \\
\end{split}
\]

\[
\begin{split}\hphantom{ \|\sum c_{a,b} Sx_a \otimes Ty_b\|^p}
&\sim \sum_{q\in Q} \int_0^1 \max_{(P_1,(w_a''))\in \mathcal P} \left( \sum_{p\in P_1}
\left(\sum_{a\in p}(\sum_{b\in q}  c_{a,b}w'_b
r_b(u))^2{w_a''}^2\right)^{p/2}\right)\,du \\
\intertext{(since $S$ is an isomorphism)}\\
&\ge \sum_{q\in Q} \int_0^1 \left( \sum_{p\in P}
\left(\sum_{a\in p}(\sum_{b\in q}  c_{a,b}w'_b
r_b(u))^2w_a^2\right)^{p/2}\right)\,du \\
&= \sum_{q\in Q} \sum_{p\in P}\int_0^1\left(\sum_{a\in
p}(\sum_{b\in q}  c_{a,b}w'_b w_a 
r_b(u))^2\right)^{p/2}\,du \\
&\sim  \sum_{q\in Q} \sum_{p\in
P}\int_0^1\int_0^1\left|\sum_{b\in q}\sum_{a\in 
p}c_{a,b}w'_b w_a r_b(u)r'_a(w)\right |^p\,dw\,du\\
\intertext{(where $(r'_a)_{a\in A}$ is another  sequence of Rademacher
functions, independent of $(r_b)_{b\in B}$)}\\
&\sim  \sum_{q\in Q} \sum_{p\in
P}\left(\sum_{b\in q}\sum_{a\in  
p}c_{a,b}^2(w'_b)^2 w_a^2\right)^{p/2}
\end{split}
\]
\end{proof}

In the proof above there are two inequalities which result from taking only
one of the partitions and weights defining the norm. If the norm is
equivalent to a norm given by finitely many partitions and weights, then we
can remove the inequality lines, insert maximums,
and complete the argument as above except that at a few places we must 
interchange the maximum with an integral. This is possible with a constant
depending on the number of partitions and weights. Therefore we have

\begin{corollary} If $X$ and $Y$ have norms given by finitely many
partitions and weights, then the $L_p$ tensor product of $X$ and $Y$ also
has a norm given by finitely many partitions and weights.
\end{corollary}

\begin{corollary} With the same hypothesis as in Proposition \ref{tensor},
 there is a constant $C'>0$ such that
for any partition $R$ of $A\times B$ and partitions and weights
$(P_r,(w_{r,a}))_{r\in R}$ and $(Q_r,(w'_{r,a}))_{r\in R}$ such that
for each $r \in R$, there exist $p_r\in P_r$ and $q_r \in Q_r$ with $r
\subset p_r \times q_r$,
$$ \|\sum c_{a,b} Sx_a \otimes Ty_b\| \ge C' \left(\sum_{r\in R}
\left( \sum_{a,b \in r} |c_{a,b}|^2 (w_{r,a})^2 (w_{r,b}')^2\right
)^{p/2}\right
)^{1/p}.$$
\end{corollary}

\begin{proof} (Sketch)
Because $SX \otimes TY$ is a subspace of $L_p$, this follows by the argument
given in the proof of Lemma
\ref{lower estimate}.
\end{proof}
\section{An Isomorph of $\bf{\otimes_{j=1}^n X_{p}}$}
\label{Distance}
In this section, we construct an example which demonstrates the
difference between a norm given by partitions and weights and the
 corresponding envelope norm. We also obtain an estimate of the distance
between a certain Banach space $Y_n$, isomorphic to $\otimes_{j=1}^n
X_{p} $, with norm given by partitions and weights,
 and any subspace of $L_p$. 
Finally we give an example of a Banach space with norm given by
partitions and weights which is 
not isomorphic to a subspace of $L_p$.\\

We will define $Y_n$ to be a Banach space with norm given by partitions and
weights which has essentially the same form as the norm on the
sequence space realization of $\otimes_{j=1}^n X_{p}$ introduced by Schechtman 
in 1975, {\bf{[S]}}. First we will estimate the distance between $Y_n$
and $Y_n$ with the associated envelope norm for the case $n=3$. Then 
for any $n\in \mathbb{N}$
we can easily extend the argument to $Y_n$ with the original norm given by 
partitions and weights and $Y_n$ with the  
corresponding envelope norm. Consequently 
we prove that not every sequence space with norm given by partitions 
and weights is isomorphic to a subspace of $L_p$ and the envelope norm 
on the sequence space realization of $\otimes_{j=1}^n X_{p}$ may be a better
choice for some purposes.  \\

 We will define $Y_3$ on
 ${\mathbb{N}}^2 \times {\mathbb{N}}^2\times {\mathbb{N}}^2$.
 Let $(w_i)_{i=1}^{\infty}$ be a sequence of weights such that 
$w_i \rightarrow 0$ as $i\rightarrow \infty$.
 Let $\bf{i_1}, \bf{i_2}$,and $\bf{i_3}$ represent indices for the first,
 second and third pair of coordinates, respectively.
Define weights on ${\mathbb{N}}^2$ by $w_{\bf{i}}=w_{(m,n)}=w_m$ where 
${\bf{i}}=(m,n)$ for all $m,n\in \mathbb{N}$.
 Let $(e_{\bf{i_1,i_2,i_3}})_{{\bf{i_1,i_2,i_3}}\in {\mathbb{N}}^2}$ be the
natural unit vector basis of $Y_3$.
The partitions of 
${\mathbb{N}}^2 \times {\mathbb{N}}^2\times {\mathbb{N}}^2$ 
and corresponding weights are given as follows,
\begin{alignat*}{3}
P_0&=\{{\mathbb{N}}^2\times{\mathbb{N}}^2\times
{\mathbb{N}}^2\}
&W_0&=(w_{\bf{i_1}} \cdot w_{\bf{i_2}} \cdot w_{\bf{i_3}})\\
P_1&=\{\{(m,n)\}\times{\mathbb{N}}^2\times{\mathbb{N}}^2: m,n\in{\mathbb{N}}
\}
&W_1&=(1\cdot w_{\bf{i_2}} \cdot w_{\bf{i_3}})\\
P_2&=\{{\mathbb{N}}^2\times\{(n,m)\}\times{\mathbb{N}}^2:
m,n\in{\mathbb{N}}\}
&W_2&=(w_{\bf{i_1}} \cdot 1 \cdot w_{\bf{i_3}})\\
P_3&=\{{\mathbb{N}}^2\times{\mathbb{N}}^2\times\{(m,n)\}:
m,n\in{\mathbb{N}}\}
&W_3&=(w_{\bf{i_1}} \cdot w_{\bf{i_2}} \cdot 1)\\
P_4&=\{\{(m,n)\}\times\{(s,t)\}\times{\mathbb{N}}^2:
m,n,s,t\in{\mathbb{N}}\}
&W_4&=(1\cdot 1 \cdot w_{\bf{i_3}})\\
P_5&=\{{\mathbb{N}}^2\times\{(m,n)\}\times
\{(s,t)\}:m,n,s,t\in{\mathbb{N}}\}
&W_5&=(w_{\bf{i_1}} \cdot 1 \cdot 1)\\
P_6&=\{\{(m,n)\}\times{\mathbb{N}}^2\times\{(s,t)\}:
 m,n,s,t \in {\mathbb{N}}\}\phantom{xxx}
&W_6&=(1\cdot w_{\bf{i_2}} \cdot 1)\\
P_7&=\{\{(l,m,n,s,t,u)\}\;\text{for}\;
l,m,n,s,t,u\in{\mathbb{N}}\}
&W_7&=(1\cdot 1 \cdot 1)
\end{alignat*}
Then the norm on $Y_3$ can be calculated by\\
\begin{equation*}\left\|\sum_{\bf{i_1,i_2,i_3}}
a_{\bf{i_1,i_2,i_3}}e_{\bf{i_1,i_2,i_3}}\right\|_{Y_3}
=\max_{I\subset \{1,2,3\}}\left\{
\left(\sum_{{\bf{i_k}}:k\in I}
\left(\sum_{{\bf{i_l}}:l\in {I^{c}}}|a_{\bf{i_1,i_2,i_3}}|^2
\prod_{l\in{I^{c}}}(w_{\bf{i_l}})^{2}\right)^{\frac{p}{2}}
\right)^{\frac{1}{p}}
\right\}
\end{equation*}
\begin{align*}
=\max\left\{\right . 
&\left(\sum_{\bf{i_1,i_2,i_3}}|a_{\bf{i_1,i_2,i_3}}|^2
(w_{\bf{i_1}})^{2}(w_{\bf{i_2}})^{2}(w_{\bf{i_3}})^{2}
\right)^{\frac{1}{2}},
\left(\sum_{\bf{i_1}}\left(\sum_{\bf{i_2,i_3}}
|a_{\bf{i_1,i_2,i_3}}|^2
(w_{\bf{i_2}})^{2}(w_{\bf{i_3}})^{2}\right)^{\frac{p}{2}}
\right)^{\frac{1}{p}},\\
&\left(\sum_{\bf{i_2}}\left(\sum_{\bf{i_1,i_3}}
|a_{\bf{i_1,i_2,i_3}}|^2
(w_{\bf{i_1}})^{2}(w_{\bf{i_3}})^{2}\right)^{\frac{p}{2}}
\right)^{\frac{1}{p}},
\left(\sum_{\bf{i_3}}\left(\sum_{\bf{i_1,i_2}}
|a_{\bf{i_1,i_2,i_3}}|^2
(w_{\bf{i_1}})^{2}(w_{\bf{i_2}})^{2}\right)^{\frac{p}{2}}
\right)^{\frac{1}{p}},\\
&\left(\sum_{\bf{i_1,i_2}}\left(\sum_{\bf{i_3}}
|a_{\bf{i_1,i_2,i_3}}|^2
(w_{\bf{i_3}})^{2}\right)^{\frac{p}{2}}
\right)^{\frac{1}{p}},
\left(\sum_{\bf{i_2,i_3}}\left(\sum_{\bf{i_1}}
|a_{\bf{i_1,i_2,i_3}}|^2
(w_{\bf{i_1}})^{2}\right)^{\frac{p}{2}}
\right)^{\frac{1}{p}},\\
&\left .\left(\sum_{\bf{i_1,i_3}}\left(\sum_{\bf{i_2}}
|a_{\bf{i_1,i_2,i_3}}|^2
(w_{\bf{i_2}})^{2}\right)^{\frac{p}{2}}
\right)^{\frac{1}{p}},
\left(\sum_{\bf{i_1,i_2,i_3}}
|a_{\bf{i_1,i_2,i_3}}|^p
\right)^{\frac{1}{p}}
\right\}
=\max\{S_i\}_{i=0}^{7}
\end{align*}
where $S_i$ for $i=0,1,\ldots,7$, are the sums
in the previous expression in the same order.

Let $Z_3$ denote $Y_3$ with the envelope norm generated.
Next we will show that 
$$
\sup_{x\in Y_3} 
\frac{\|x\|_{Z_3}}{\|x\|_{Y_3}} \ge 3^{\frac{1}{p}}.
$$

Since $w_i \rightarrow 0$ as $i \rightarrow \infty$,
 then for any given $\epsilon $, $0<\epsilon \leq 3$,
 there exists an $N$, such that if $n > N$, 
$w_n<\left(\frac{\epsilon}{3}\right)^{\frac{1}{2}}\leq
\left(\frac{\epsilon}{3}\right)^{\frac{1}{p}}$.
Let $n_1, n_2, n_3 > N$.
Choose integers $K_1, K_2, K_3$, such that
$$w_{n_1}K_{1}^{\frac{1}{2}-\frac{1}{p}}
> \left(\frac{3}{\epsilon}\right)^{\frac{1}{p}}\ge 1, \quad
w_{n_2}K_{2}^{\frac{1}{2}-\frac{1}{p}}
> \left(\frac{3}{\epsilon}\right)^{\frac{1}{p}}\ge 1,\quad
w_{n_3}K_{3}^{\frac{1}{2}-\frac{1}{p}}
> \left(\frac{3}{\epsilon}\right)^{\frac{1}{p}}\ge 1.$$ 
Now take three blocks with constant coefficients as follows:\newline
a block of size $K_1$, $x_1$, with coefficient
$(w_{n_1})^{-1}K_{1}^{-\frac{1}{2}}$ and support
\begin{alignat*}{2}
&\{(n_1,1, n_2,1, n_3,1),
(n_1,2, n_2,1,n_3,1),
&&\dots,
(n_1, K_1, n_2,1, n_3,1)\},\\
\intertext{a block of size $ K_2$, $x_2$, with coefficient
$(w_{n_2})^{-1} K_{2}^{-\frac{1}{2}}$ and support}\\
&\{(n_1,K_1+1, n_2,2, n_3,2),
(n_1,K_1+1,&&n_2,3, n_3,2),\\
&&&\dots,
(n_1, K_1+1, n_2,K_2+1, n_3,2)\},\\
\intertext{and a block of size $K_3$, $x_3$, with coefficient
$(w_{n_3})^{-1} K_{3}^{-\frac{1}{2}}$ and support}\\
&\{(n_1, K_1+2, n_2,K_2+2, n_3,3),
(n_1,K_1&&+2,n_2,K_2+2, n_3,4),\\
&&&\dots,
(n_1, K_1+2, n_2,K_2+2, n_3,K_3+2)\}.
\end{alignat*}

Now we estimate the eight sums to get an estimate of the norm of
 the element
\begin{multline}\label{vector}
x_1+x_2+x_3 =\sum_{{\bf{i_1}}=(n_1,1)}^{(n_1,K_1)}
 w_{n_1}^{-1}K_{1}^{-\frac{1}{2}}
 e_{{\bf{i_1}},n_2,1,n_3,1} \\ 
+\sum_{{\bf{i_2}}=(n_2,2)}^{(n_2,K_2+1)} 
w_{n_2}^{-1}K_{2}^{-\frac{1}{2}}
 e_{n_1,K_1+1,{\bf{i_2}},n_3,2}
+\sum_{{\bf{i_3}}=(n_3,3)}^{(n_3,K_3+2)}
  w_{n_3}^{-1}K_{3}^{-\frac{1}{2}}
 e_{n_1,K_1+2,n_2,K_2+2,{\bf{i_3}}}.
\end{multline}
\begin{align*}
S_0&=\left[(w_{n_1})^{-2}K_{1}^{-1}
(w_{n_1})^{2}
(w_{n_2})^{2}
(w_{n_3})^{2}K_1\right .
 +(w_{n_2})^{-2}K_{2}^{-1}
(w_{n_1})^{2}(w_{n_2})^{2}(w_{n_3})^{2}K_2\\
&\left .+(w_{n_3})^{-2}K_{3}^{-1}
(w_{n_1})^{2}(w_{n_2})^{2}(w_{n_3})^{2}K_3
\right]^{\frac{1}{2}}\\
& =\left[(w_{n_2})^{2}(w_{n_3})^{2}
+(w_{n_1})^{2}(w_{n_3})^{2}
+(w_{n_1})^{2}(w_{n_2})^{2}
\right]^{\frac{1}{2}}
 <{\epsilon}^{\frac{1}{2}}.\\
S_1&=[((w_{n_1})^{-2}K_{1}^{-1}
(w_{n_2})^{2}(w_{n_3})^{2})^{\frac{p}{2}}K_1
+((w_{n_2})^{-2}K_{2}^{-1}
(w_{n_2})^{2}(w_{n_3})^{2}K_2)^{\frac{p}{2}}\\
&+((w_{n_3})^{-2}K_{3}^{-1}
(w_{n_2})^{2}(w_{n_3})^{2}K_3)^{\frac{p}{2}}
]^{\frac{1}{p}}\\
&=[(w_{n_1}K_{1}^{\frac{1}{2}-\frac{1}{p}})^{-p}
(w_{n_2})^{p}(w_{n_3})^{p}+(w_{n_3})^{p}
+(w_{n_2})^{p}]^{\frac{1}{p}}<{\epsilon}^{\frac{1}{p}}.
\end{align*}
Similarly we have 
$S_2<{\epsilon}^{\frac{1}{p}}$
and
$S_3<{\epsilon}^{\frac{1}{p}}$.\\
\begin{align*}
S_4&=\big(K_1
((w_{n_1})^{-2}K_{1}^{-1}
(w_{n_3})^2)^{\frac{p}{2}}
+K_2
((w_{n_2})^{-2}K_{2}^{-1}
(w_{n_3})^2)^{\frac{p}{2}}
+((w_{n_3})^{-2}K_{3}^{-1}
(w_{n_3})^2K_3)^{\frac{p}{2}}
\big)^{\frac{1}{p}}\\
&=\big((w_{n_1}
K_{1}^{\frac{1}{2}-\frac{1}{p}})^{-p}
(w_{n_3})^p
+(w_{n_2}
K_{2}^{\frac{1}{2}-\frac{1}{p}})^{-p}
(w_{n_3})^p
+1\big)^{\frac{1}{p}}
<(\epsilon + 1)^{\frac 1 p}.
\end{align*}
Similarly we have 
$S_5<(\epsilon + 1)^{\frac 1 p}$
and
$S_6<(\epsilon + 1)^{\frac 1 p}$.\\
\begin{align*}
&S_7=((w_{n_1}K_{1}^{\frac{1}{2}-\frac{1}{p}})^{-p}
+(w_{n_2}K_{2}^{\frac{1}{2}-\frac{1}{p}})^{-p}
+(w_{n_3}K_{3}^{\frac{1}{2}-\frac{1}{p}})^{-p})^{\frac{1}{p}}
< {\epsilon}^{\frac 1 p}.
\end{align*}
Since $\epsilon $ can be arbitrary small, if we take maximum of
 these eight sums, the norm will be as close to 1 as we want.\\

Now let us look at the envelope norm of the element $x_1+x_2+x_3$ in (\ref{vector}).\\ 

Let $Q$ be a partition of 
 ${\mathbb{N}}^2 \times {\mathbb{N}}^2\times {\mathbb{N}}^2$
 such that the support of each of the above three blocks is an 
element of $Q$. (The other sets in the partition do not matter.)
Let $\mathcal{P}$ be the given family of weights and partitions, i.e., 
${\mathcal{P}}=\{ (P_i,W_i): i=0,1,...,7\}$. 
Let $T: Q\rightarrow \mathcal{P}$ be a map such that
$T(\text{supp }x_1)=(P_6,W_6)$, 
$T(\text{supp }x_2)=(P_5,W_5)$, and
$T(\text{supp }x_3)=(P_4,W_4)$. 
Then the envelope norm of $x_1+x_2+x_3$ can be 
estimated from below using $P(Q, T)$ 
\begin{align*}
&\left\|\left|\sum_{{\bf{i_1}}=(n_1,1)}^{(n_1,K_1)}
 w_{n_1}^{-1}K_{1}^{-\frac{1}{2}}
 e_{{\bf{i_1}},n_2,1,n_3,1}
+\sum_{{\bf{i_2}}=(n_2,2)}^{(n_2,K_2+1)} 
w_{n_2}^{-1}K_{2}^{-\frac{1}{2}}
 e_{n_1,K_1+1,{\bf{i_2}},n_3,2}\right .\right .\\
&\left .\left.+\sum_{{\bf{i_3}}=(n_3,3)}^{(n_3,K_3+2)}
  w_{n_3}^{-1}K_{3}^{-\frac{1}{2}}
 e_{n_1,K_1+2,n_2,K_2+2,{\bf{i_3}}}\right\|\right|
\ge\left( \left(\sum_{i_1=(n_1,1)}^{(n_1,
 K_1)}(w_{n_1}^{-1}K_{1}^{-\frac{1}{2}})^2
w_{i_1}^{2}\right)^{\frac p 2}\right .\\
&+\left(\sum_{i_2=(n_2,2)}^{(n_2,
 K_2+1)}(w_{n_2}^{-1}K_{2}^{-\frac{1}{2}})^2
w_{i_2}^{2}\right)^{\frac p 2}
+\left.\left(\sum_{i_3=(n_3,3)}^{(n_3,
 K_3+2)}(w_{n_3}^{-1}K_{3}^{-\frac{1}{2}})^2
w_{i_3}^{2}\right)^{\frac p 2}\right)^{\frac 1 p}\\
&\ge
\left(\left((w_{n_1})^{-2}K_{1}^{-1}(w_{n_1})^{2}K_1\right)^{\frac p
    2}
+\left((w_{n_2})^{-2}K_{2}^{-1}(w_{n_2})^{2}K_2\right)^{\frac p
 2}\right .
+\left .\left((w_{n_3})^{-2}K_{3}^{-1}(w_{n_3})^{2}K_3\right)^{\frac
 p 2}
\right)^{\frac 1 p}\\
&= 3^{\frac{1}{p}}
\end{align*}
Hence the envelope norm on $Y_3$ is at best $3^{\frac 1 p}$ 
equivalent to the given norm.\\

Next we will describe how this computation can be generalized. We define
  $Y_n$ for any $n\in \mathbb{N}$  on 
${\mathbb{N}}^2\times {\mathbb{N}}^2\times \hdots \times
{\mathbb{N}}^2$(n times).
Let $w_{\bf{i}}= w_{s,t}=w_s$ for ${\bf{i}}=(s,t),\; s,t\in
  \mathbb{N}$ as above.
Let $I\subset \{1,2,\hdots,n\}$. Define\\
 $$\begin{displaystyle}P_I=\left\{
\prod_{k=1}^{n} A_k:\; \text{where}\; A_k={\mathbb{N}}^2,
k\notin I;\;A_k=\{(m_k, l_k)\}, k\in I, \; m_k, l_k \in
\mathbb{N}\right \}\end{displaystyle}$$
 and
$$\begin{displaystyle} W_I=\left(\prod_{k \notin I}w_{\bf{i_k}}
\right).\end{displaystyle}$$ 
For a given sequence $(w_i)$ such that
 $w_i \rightarrow 0$ as $i \rightarrow \infty$ and 
any $0 <\epsilon \le 1$, there exists $N$, such that
 if $m> N$, then 
$w_m<
\left(\frac{\epsilon}{n}\right)^{\frac 1 2}\le
\left(\frac{\epsilon}{n}\right)^{\frac 1 p}.
$ 

Let $m_1,\ldots, m_n > N$.
We choose $n$ blocks with size $K_l$ for $l=1,\ldots, n$ 
in $({\mathbb{N}}^2)^n$ so that 
$w_{m_l}{K_l}^{{\frac 1 2}-{\frac 1 p}}> (\frac n \epsilon )^{\frac 1 p}$.\\

The block of size $K_{l+1}$ would have coefficient
$w_{m_{l+1}}^{-1}K_{l+1}^{-\frac 1
2}$ and have support
\begin{align*}
&(m_1,K_1+l,m_2,K_2+l,\ldots, m_{l+1},l+1, \ldots,m_n, l+1)\\
&(m_1,K_1+l,m_2,K_2+l,\ldots, m_{l+1},l+2, \ldots,m_n, l+1)\\
&\phantom{xxxxxxxxxxxxxxxx}\vdots\\
&(m_1,K_1+l,m_2,K_2+l,\ldots, m_{l+1},l+K_{l+1}, \ldots,m_n, l+1)
\end{align*}
where $0\le l \le n-1$.\\

By applying similar arguments to that for $n=3$ 
we have that the value of 
 the envelope norm of the sum of these blocks is at least
$n^{\frac{1}{p}}$ while the value of the norm given 
by partitions and weights remains approximately 1.\\

\begin{theorem}
The distance from $Y_n$ to a subspace of $L_p$ goes to $\infty$ with
$n$, i.e., there is a sequence $(K(n))$, $K(n)\rightarrow \infty$,
such that for all isomorphisms $T: \; Y_n \rightarrow Z\subset L_p$, 
$\|T\|\|T^{-1}\|\geq K(n)$.
\end{theorem}

Proof: If  $T: \; Y_n \rightarrow Z\subset L_p$ is an isomorphism, then 
 by Theorem \ref{equivalent norm}, the norm of $Y_n$ 
given by partitions and weights is equivalent to the envelope norm
with a constant depending on $\|T\|\|T^{-1}\|$. 
Since the envelope norm of some element of $Y_n$ of norm 1 
has value at least $n^{\frac 1 p}$,
 then $\begin{displaystyle}\|T\|\|T^{-1}\|\geq {\lambda}^{-1}n^{\frac 1 p}
\geq \frac{n^{\frac{1}{p}}}{\|T\|\|T^{-1}\|}\end{displaystyle}$\qed

\begin{corollary} The distance between $Y_n$ and $\otimes_{j=1}^n X_{p}$ 
goes to $\infty$ with $n$. 
\end{corollary}
\begin{corollary}
$(\sum_n Y_n)_{\ell_p}$ with norm given by partitions and weights is not
isomorphic to a subspace of $L_p$.
\end{corollary}
\section{Remarks and Open Problems}

We introduced the envelope property to show that there is a space with
norm given by partitions and weights which is not isomorphic to subspace
of $L_p$. It seems unlikely that this property alone determines whether
a space with
norm given by partitions and weights  is isomorphic to a subspace of $L_p$.

\begin{question} Is there a space with
norm given by partitions and weights  which has the envelope property
but is not isomorphic to a subspace of $L_p$?
\end{question}

\begin{question} What are necessary and sufficient conditions for a space
with
norm given by partitions and weights  to be isomorphic to a subspace of
$L_p$?
\end{question}

We showed that the tensor product of two spaces with norms given by
finitely many partitions and weights is isomorphic to a space with
norm given by partitions and weights.

\begin{question} Suppose $X$ and $Y$ have norms given by partitions and
weights  and are each isomorphic to a subspace of $L_p$. Is $X \otimes Y$
isomorphic to a subspace of $L_p$? What if $X$ and $Y$ are each isomorphic
to a complemented subspace of $L_p$?
\end{question}

The construction of uncountably many complemented subspaces of $L_p$ given
 by Bourgain, Rosenthal and Schechtman is based on two fundamental
operations. If $X$ and $Y$ are subspaces of $L_p[0,1]$ then a
distributional version of the $\ell_p$-sum is used, $X\oplus_p Y,$ with
each space being isometrically mapped to a space supported on half the
interval in a canonical way. The second operation
is used with sequences of subspaces of $L_p$. Let $(X_n)$ be such a
sequence such that each $X_n$ contains the constant functions and let $X_{n,0}$
denote the subspace of mean zero functions in $X_n$. An infinite sum is
created by placing (isometrically)
$X_{n,0}$ onto an infinite product of probability
spaces as functions depending only on the $n$th coordinate.
The space $(\sum X_n)_I$ is the span of these transported copies of
$X_n$, $n=1,2,\dots$ and the constant functions. 

With these two operations an induction on $\omega_1$ is used to
construct the spaces. Let $R_p^0$ be the constant functions on $[0,1]$.
If $R_p^\alpha$ has
been defined, define $R_p^{\alpha+1}=R_p^\alpha \oplus_p R_p^\alpha.$
If $\alpha$ is a limit ordinal, define $R_p^\alpha =(\sum_{\beta<\alpha}
R_p^\beta)_I.$

These operations and the construction are investigated in {\bf[A]}.
In particular it is shown that there that the constant functions do not play an
important role in the construction. If we make a few adjustments we can
mimic this construction using $(p,2,W)$ sums. For the $\ell_p$-sum of
two spaces, we use the $(p,2,W_1)$ sum with two equal weights,
$W_1=(2^{(2-p)/(2p)},2^{(2-p)/(2p)})$. For the independent sum we use the
$(p,2,W_2)$ sum with $W_2$ equal to the constantly $1$ sequence. In each case
we will also assume that we take the envelope norm generated. In addition
 in order to strengthen the correspondence between this construction
and very distributional construction of $R_p^\alpha$,
we will use at each step the fact that the space $R_{p,0}^\alpha$ has an
unconditional  basis which is orthogonal in $L_2$.
Define $Y_p^0 =[1_{[0,1/2)}-1_{[1/2,1]}].$  
This is just  a one dimensional space. Let $X_p^0$ be the
sequences of length one.
If $Y_p^\alpha$ and $X_p^\alpha$ have been defined, define
$Y_p^{\alpha+1}=Y_p^\alpha \oplus_p Y_p^\alpha$ and
$X_p^{\alpha+1}=(X_p^\alpha,X_p^\alpha)_{(p,2,W_1)}.$
If $\alpha$ is a limit ordinal, let $Y_p^\alpha =(\sum_{\beta<\alpha}
Y_p^\beta)_I$ and $X_p^{\alpha+1}=\sum_{\beta<\alpha}
X_p^\beta)_{(p,2,W_2)}$.

It follows from the results in Chapter 2 of {\bf[A]} that for 
$\alpha<\omega^2$
that $R_p^\alpha$, $Y_p^\alpha$, and $X_p^{\alpha}$ are isomorphic. The
choice of weights $W_1$ and $W_2$ is such that the natural mapping from
$Y_p^\alpha,\|\cdot\|_2$ to $X_p^\alpha,\|\cdot\|_2$ will be an isometry.

\begin{question}
Is $X_p^\alpha$ isomorphic to $R_p^\alpha$ for all $\alpha<\omega_1$?
\end{question}

We think it unlikely to be true
 but it would be nice to know the answer to the
following simple question.

\begin{question} Is $L_p$ isomorphic to a space with norm given by
partitions and weights?
\end{question}

Finally we note that the envelope norm suggests that a useful alternate form of
Rosenthal's inequality for a sequence of mean zero independent random
variables $(f_n)$ might be
\begin{multline*}c_p \max_{Q\subset \mathbb N}\{ \left(\sum_{n\notin Q}
\|f_n\|^p_p+\left (\sum_{n\in Q}
\|f_n\|^2_2 \right )^{p/2}\right )^{1/p}\} \le \|\sum_n f_n \|_p \\
\le C_p  \max_{Q\subset \mathbb N}\{ \left (\sum_{n\notin Q} 
\|f_n\|^p_p+\left (\sum_{n\in Q}
\|f_n\|^2_2\right)^{p/2}\right )^{1/p}\}
 .\end{multline*}

\begin{question} What are the best constants in this form of Rosenthal's
inequality?
\end{question}

(See {\bf [J-S-Z]} for results on the constants in the original form.)
\vskip 1.5 cm

{\center{{\bf{References}}}}\\
{\bf{[A]}}\; D. Alspach, \textit{Tensor products and independent sums of 
${\mathcal{L}}_p$-spaces, $1<p<\infty$},
Mem. Amer. Math. Soc. {\bf{138}} (660) (1999), viii+77.\\
{\bf{[A-O]}}\; D. Alspach and E. Odell, \textit{$L_p$ Spaces, Handbook of the
Geometry of Banach Spaces}, Vol.1, W.B.Johnson and J. Lindenstrauss,
eds, Elsevier, Amsterdam (2001), 123-159.\\
{\bf{[B-R-S]}}\; J. Bourgain, H.P. Rosenthal and G. Schechtman,\textit{ An
ordinal $L_p$-index for Banach spaces, with application to
complemented subspaces of $L_p$}, Ann. of Math. (2) {\bf{114}}(2)
(1981), 193-228.\\
{\bf{[F]}}\; G. Force, \textit{Constructions of ${\mathcal{L}}_p$-spaces,
$1<p\ne 2< \infty$}, Dissertation, Oklahoma State University,
Stillwater, Oklahoma, 1995, Available from ArXiv.org as math.FA/9512208\\
{\bf{[J-L]}}\;\textit{Handbook of the geometry of Banach Spaces},
 Vol.1, W.B.Johnson and J. Lindenstrauss,
eds, Elsevier, Amsterdam (2001).\\
{\bf{[J-S-Z]}}\; W. B. Johnson, G. Schechtman, and J. Zinn,
  \textit{Best constants in moment inequalities for linear
  combinations of independent and exchangeable random variables,}
  Ann. Probab. {\bf{13}} (1985), 234-253.\\
{\bf{[L-P]}}\; J. Lindenstraus and A. Pelczynski,\textit{ Absolutely summing
operators in ${\mathcal{L}}_p$-spaces and their applications}, Studia
Math. {\bf{29}} (1968), 275-326.\\
{\bf{[L-R]}}\;  J. Lindenstraus and  H.P. Rosenthal, \textit{The
${\mathcal{L}}_p$ Spaces}, Israel J. Math. {\bf{7}} (1969), 325-349.\\
{\bf{[L-T-1]}}\; J. Lindenstrauss and L. Tzafriri,\textit{ Classical Banach
Spaces. I, Sequences Spaces}, Ergebnisse der Mathematik und ihrer 
Grenzgebiete 92, Springer-Verlag, Berlin (1977).\\
{\bf{[L-T-2]}}\; J. Lindenstrauss and L. Tzafriri,\textit{ Classical Banach
Spaces. II }, Ergebnisse der Mathematik und ihrer 
Grenzgebiete 97, Springer-Verlag, Berlin (1979).\\
{\bf{[P]}}\; Pelczynski, \textit{Projections in certain Banach
  spaces}, Studia Math. {\bf{19}} (1960) 209-228 \\
{\bf{[R]}}\; H.P. Rosenthal,\textit{ On the subspaces of $L_p(p>2)$ spanned by 
sequences of independent random variables}, Israel J. Math. {\bf{8}}
(1970), 273-303.\\
{\bf{[S]}}\; G. Schechtman,\textit{ Examples of ${\mathcal{L}}_p$ spaces 
($1<p\ne 2 <\infty$)}, Israel J. Math. {\bf{22}} (1) (1975),
  138-147.

\end{document}